\documentclass[leqno]{article}
\usepackage{amsmath}
\usepackage{amsfonts}
\usepackage{amssymb}

\usepackage[top=30truemm,bottom=30truemm,left=25truemm,right=25truemm]{geometry}

\newcommand{\qed}{\hbox{\rule[-2pt]{3pt}{6pt}}}

\def\Epsilon{{\rm E}}

\renewcommand{\Im}{\operatorname{Im}}
\renewcommand{\Re}{\operatorname{Re}}
\renewcommand{\det}{\operatorname{det}}
\renewcommand{\mod}{\operatorname{mod}}

\begin{document}

\hfil {\bf Pullback formula for vector valued Siegel modular forms and its applications} \hfil

\bigskip

\hfil {Noritomo Kozima} \hfil

\bigskip

\bigskip
\noindent{\bf 1. Introduction} \par

 Let $E^n_k(Z,s)$ be the Eisenstein series of degree $n$ and weight $k$ with a complex variable $s$, i.e.,
\[ E^n_k(Z,s) := \det(\Im(Z))^s \sum_{(C,D)} \det(CZ+D)^{-k} \left|\det(CZ+D)\right|^{-2s} \]
where $(C,D)$ runs over coprime symmetric pairs of degree $n$. The Eisenstein series converges for $k+2\Re(s)>n+1$. As is well known, $E^n_k(Z,s)$ has meromorphic continuation to the whole $s$-plane and satisfies a functional equation.

 On the other hand, let $f$ be a Siegel cuspform of weight $k$ with respect to $Sp(q,\mathbb{Z})$ (size $2q$). Suppose $f$ is an eigenform, i.e., a non-zero common eigenfunction of the Hecke algebra. Then taking Petersson inner product of $f(W)$ and $E^{p+q}_k\Bigl(\begin{pmatrix} -\overline{Z}^{(p)} & 0 \\ 0 & W^{(q)}\end{pmatrix},\overline{s}\Bigr)$ in variable $W$, B\"ocherer [1] proved the following:
\[ \left(f,E^{p+q}_k\Bigl(\begin{pmatrix} -\overline{Z} & 0 \\ 0 & * \end{pmatrix},\overline{s}\Bigr)\right) = c_q(s)\, L(2s+k-q,f,\underline{\rm St})\, [f]^p_q(Z,s), \]
where $c_q(s)$ is a product of the gamma functions, $L(s,f,\underline{\rm St})$ the standard $L$-function attached to $f$, and $[f]^p_q(Z,s)$ the Klingen Eisenstein series of $f$.

Using the above identity, he investigated analytic continuation and functional equation with respect to the standard $L$-functions and Klingen Eisenstein series.

This theory was generalized in symmetric tensor valued case ([3], [17], [19], [11]) and alternating tensor valued case ([18], [12]). In this paper, under a certain condition, we generalize this theory in any vector valued cases. The main theorem is as follows.

\bigskip
{\bf Theorem 5.2.}
{\sl
 Let\/ $k$ be even and\/ $k+2\Re(s)>p+q+1$. If\/ $V_p^*\otimes V_{q*}$-valued polynomial $\mathfrak{P}(X)$ is homogeneous in\/ $X$ with underlying set\/ $I$ and pluri-harmonic for\/ $\Delta^*(2k)$ and\/ $\Delta_*(2k)$, then for any eigenform\/ $f\in\mathfrak{S}_{\rho_q}$,
\begin{align*}
    & \left( f_*,(\mathfrak{P}(\partial)E^{p+q}_k)\Bigl(\begin{pmatrix} -\overline{Z}^{(p)} & 0 \\ 0 & * \end{pmatrix},\overline{s} \Bigr) \right)
\\  & \qquad\qquad = 2^{q(q+1-2s)+1}\, (2^{-1}i)^{\left|\rho_q\right|}\, c(s,\rho_q)\, D(k+2s,f)\, [f]^{\rho_p}_{\rho_q}(Z,s)^*.
\end{align*}
Here\/ $c(s,\rho_q)$ is a\/ $\mathbb{C}$-valued function of\/ $s$.
}
\bigskip

\bigskip
\noindent{\bf 2. Preliminaries} \par

\bigskip
\noindent{\bf 2.1. Vector valued Siegel modular forms} \par

 Let $n$ be a positive integer. Let $(\rho,V_\rho)$ be an irreducible polynomial representation of $GL(n,\mathbb{C})$ on a finite-dimensional complex vector space $V_\rho$ such that the highest weight of $\rho$ is $(\lambda_1,\lambda_2,\ldots,\lambda_n)\in\mathbb{Z}^n$ with $\lambda_1\geq\lambda_2\geq\ldots\geq\lambda_n\geq 0$. We put $\left|\rho\right|:=\sum\limits_{j=1}^n\lambda_j$. We fix a Hermitian inner product $\langle \cdot,\cdot \rangle$ on $V_\rho$ such that
\[ \langle \rho(g)v,w \rangle = \langle v,\rho({}^t\overline{g})w \rangle \quad \text{for $g\in GL(n,\mathbb{C})$, $v$, $w\in V_\rho$}. \]

 For a subring $R$ of $\mathbb{R}$, we put
\[ G^+Sp(n,R) := \left\{ g\in M(2n,R) \Bigm| {}^tg\begin{pmatrix} 0 & 1_n \\ -1_n & 0 \end{pmatrix} g = r(g)\begin{pmatrix} 0 & 1_n \\ -1_n & 0 \end{pmatrix},\quad\text{for some $r(g)>0$} \right\} \]
and
\[ Sp(n,R):=\left\{g\in G^+Sp(n,R) \mid r(g)=1 \right\}. \]
Let $\Gamma_n:=Sp(n,\mathbb{Z})$ be the Siegel modular group of degree $n$, and $\mathfrak{H}_n$ the Siegel upper half space of degree $n$. For $g=\begin{pmatrix} A^{(n)} & B^{(n)} \\ C^{(n)} & D^{(n)} \end{pmatrix}\in G^+Sp(n,\mathbb{R})$ and $Z\in\mathfrak{H}_n$, we put
\[ g\langle Z\rangle:=(AZ+B)(CZ+D)^{-1},\quad j(g,Z):=CZ+D, \]
and for a $V_\rho$-valued function $f\colon\mathfrak{H}_n\to V_\rho$,
\[ (f|_\rho g)(Z):=r(g)^{\left|\rho\right|-n(n+1)/2}\rho(j(g,Z)^{-1}) f(g\langle Z\rangle). \]
We write $|_k$ for $\rho=\det^k$

 A $C^\infty$-function $f\:\mathfrak{H}_n\to V_\rho$ is called a $V_\rho$-valued modular form of weight $\rho$ if it satisfies $f|_\rho\gamma=f$ for all $\gamma\in\Gamma_n$. The space of all such functions is denoted by $\mathfrak{M}_\rho^\infty$. The space of $V_\rho$-valued Siegel modular forms of weight $\rho$ is defined by
\[ \mathfrak{M}_\rho := \{ f\in\mathfrak{M}_\rho^\infty \mid \text{$f$ is holomorphic on $\mathfrak{H}_n$ (and its cusps)} \}, \]
and the space of cuspforms by
\[ \mathfrak{S}_\rho := \left\{ f\in\mathfrak{M}_\rho \Bigm| \lim_{\lambda\to\infty} f\bigl(\left(\begin{matrix} Z^{(n-1)} & 0 \\ 0 & i\lambda \end{matrix} \right)\bigr)=0 \quad\text{for all $Z\in\mathfrak{H}_{n-1}$} \right\}. \]
If $\rho=\det^k$, we write ${\mathfrak{M}^n_k}^\infty$, $\mathfrak{M}^n_k$, and $\mathfrak{S}^n_k$ for $\mathfrak{M}_\rho^\infty$, $\mathfrak{M}_\rho$, and $\mathfrak{S}_\rho$, respectively.

 For $f$, $g\in\mathfrak{M}_\rho^\infty$, the Petersson inner product of $f$ and $g$ is defined by
\[ (f,g):=\int_{\Gamma_n\backslash\mathfrak{H}_n} \left\langle \rho(\sqrt{\Im(Z)})f(Z), \rho(\sqrt{\Im(Z)})g(Z) \right\rangle \det(\Im(Z))^{-n-1} \,dZ \]
if the right-hand side is convergent.

\bigskip
\noindent{\bf 2.2. Hecke algebra and standard $L$-functions} \par

We consider the Hecke algebra $\mathcal{H}(\Gamma_n,G^+Sp(n,\mathbb{Q}))$. For $g\in G^+Sp(n,\mathbb{Q})$, let $\Gamma_n g\Gamma_n=\displaystyle{\bigsqcup_{j=1}^r\Gamma_n g_j}$ be a decomposition of the double coset $\Gamma_n g\Gamma_n$ into left cosets. For $f\in\mathfrak{M}_\rho$ (resp. $\mathfrak{S}_\rho$, $\mathfrak{M}_\rho^\infty$), we define the Hecke operator $(\Gamma_n g\Gamma_n)$ by
\[ f|(\Gamma_n g\Gamma_n):=\sum_{j=1}^r f|_\rho g_j. \]

 Let $f\in\mathfrak{S}_\rho$ be an eigenform. We define the standard $L$-function attached to $f$ by
\[ L(s,f,\underline{\rm St}):=\prod_p\left\{(1-p^{-s})\prod_{j=1}^n(1-\alpha_j(p)p^{-s}) (1-\alpha_j(p)^{-1}p^{-s}) \right\}^{-1}, \tag{2.1} \]
where $p$ runs over all prime numbers and $\alpha_j(p)$ ($1\leq j\leq n$) are the Satake $p$-parameters of $f$. The right-hand side of (2.1) converges absolutely and locally uniformly for $\Re(s)$ sufficiently large.

 We also define the following Dirichlet series:
\[ D(s,f):=\sum_{T\in\mathbb{T}^{(n)}} \lambda(T,f) \det(T)^{-s}, \tag{2.2} \]
where $\mathbb{T}^{(n)}$ is the set consisting of elementary divisor forms, i.e.,
\[ \font\b=cmr10 scaled \magstep2
   \def\bigzerol{\smash{\hbox{\b 0}}}
   \def\bigzerou{\smash{\lower.9ex\hbox{\b 0}}}
   \mathbb{T}^{(n)}:=\left\{ \begin{pmatrix} t_1 & & & \bigzerou \\ & t_2 & & \\ & & \ddots & \\ \bigzerol & & & t_n \end{pmatrix} \Biggm| \text{$t_1$, $t_2$, $\ldots$, $t_n$ are positive integers such that $t_1 \mid t_2 \mid \ldots \mid t_n$} \right\}, \]
and $\lambda(T,f)$ is the eigenvalue on $f$ of the Hecke operator $\Bigl(\Gamma_n \begin{pmatrix} T & 0 \\ 0 & T^{-1} \end{pmatrix} \Gamma_n \Bigr)$. Then by B\"ocherer [2], we have the following identity:
\[ \zeta(s)\prod_{j=1}^n\zeta(2s-2j)D(s,f) = L(s-n,f,\underline{\rm St}). \tag{2.3} \]

\bigskip
\noindent{\bf 2.3. Eisenstein series} \par

 For $k\in 2\mathbb{Z}_{>0}$, $s\in\mathbb{C}$ and $Z\in\mathfrak{H}_n$, we define the Eisenstein series by
\[ E^n_k(Z,s) := \det(\Im(Z))^s \sum_{\gamma\in\Gamma_{n,0}\backslash\Gamma_n} \det(j(\gamma,Z))^{-k} \left| \det(j(\gamma,Z)) \right|^{-2s}, \]
where
\[ \Gamma_{n,r}:=\left\{\begin{pmatrix} * & * \\ 0^{(n+r,n-r)} & * \end{pmatrix} \in \Gamma_n \right\}. \]
Then $E^n_k(Z,s)$ converges absolutely and locally uniformly for $k+2\Re(s)>n+1$, and $E^n_k(Z,s)\in{\mathfrak{M}^n_k}^\infty$. As is well known, $E^n_k(Z,s)$ has meromorphic continuation to the whole $s$-plane and a functional equation:

\bigskip
{\bf Theorem 2.1.} (Langlands [14], Kalinin [8], Diehl [4] and Mizumoto [15, 16])
{\sl
 Let\/ $n$ be a positive integer and\/ $k$ a positive even integer. We put
\[ \mathbb{E}^n_k(Z,\,s) := \gamma^n_k(s)\, E^n_k\left(Z,\,s-\frac{k}{2}\right), \]
where
\[ \gamma^n_k(s) := \frac{\Gamma_n\Bigl(\displaystyle{s+\frac{k}{2}}\Bigr)}{\Gamma_n(s)} \xi(2s) \prod_{j=1}^{[n/2]} \xi(4s-2j) \]
with
\[ \xi(s):=\pi^{-s/2}\Gamma\Bigl(\frac{s}{2}\Bigr)\zeta(s) \quad\text{and}\quad \displaystyle{\Gamma_n(s):=\prod_{j=1}^n \Gamma\Bigl(s-\frac{j-1}{2}\Bigr)}. \]
Then $\mathbb{E}^n_k(Z,s)$ is an entire function in $s$, and it is invariant under the substitution $s\mapsto\displaystyle{\frac{n+1}{2}-s}$.
}
\bigskip

Furthermore, it is also known that every partial derivative (in the entries of $Z$ and $\overline{Z}$) of the Eisenstein series $E^n_k(Z,s)$ is slowly increasing (locally uniformly in $s$):

\bigskip
{\bf Theorem 2.2.} (Mizumoto [16])
{\sl
 Let\/ $n$ be a positive integer and\/ $k$ a positive even integer.

\noindent {\rm (i)} For each\/ $s_0\in\mathbb{C}$, there exist constants\/ $\delta>0$ and a non-negative integer\/ $d$ depending only on\/ $n$, $k$ and\/ $s_0$, such that
\[ (s-s_0)^d\,E^n_k(Z,s) \]
is holomorphic in\/ $s$ for\/ $\left|s-s_0\right|<\delta$, and is\/ $C^\infty$ in\/ $(\Re(Z),\Im(Z))$.

\noindent {\rm (ii)} Furthermore, for given\/ $\varepsilon>0$ and non-negative integer $N$, there exist constants\/ $\alpha>0$ and\/ $\beta>0$ depending only on\/ $n$, $k$, $d$, $s_0$, $\varepsilon$, $\delta$ and\/ $N$ such that
\[ \left|(s-s_0)^d\,\mathcal{D}E^n_k(Z,s)\right| \leq \alpha\det(\Im(Z))^\beta \]
for\/ $Y\geq\varepsilon 1_n$ and\/ $\left|s-s_0\right|<\delta$, where\/ $\mathcal{D}$ is an arbitary monomial of degree\/ $N$ in $\displaystyle{\frac{\partial}{\partial z_{\mu\nu}}}$ and\/ $\displaystyle{\frac{\partial}{\partial \overline{z}_{\mu\nu}}}$ ($1\leq\mu,\nu\leq n$). Here\/ $z_{\mu\nu}$ denotes the\/ $(\mu,\nu)$-entry of\/ $Z$.
}
\bigskip

\noindent Remark. The assertion above for the case $N=0$ has been proved by Langlands [14] and Kalinin [8].

\bigskip

{\bf Theorem 2.3.} (Feit [5])
{\sl
 Let\/ $n$ be a positive integer and\/ $k$ a positive even integer. We put
\[ \tilde{E}^n_k(Z,s) := \zeta(k+2s)\, \prod_{j=1}^{[n/2]} \zeta(2k+4s-2j)\, E^n_k(Z,s). \]

\noindent{\rm (i)} If\/ $k\geq n/2$ and\/ $n\not\equiv 0$ $(\mod 4)$, then $\tilde{E}^n_k(Z,s)$ is holomorphic for $k+2\Re(s)>[n/2]$.

\noindent{\rm (ii)} If\/ $k\geq n/2$ and\/ $n\equiv 0$ $(\mod 4)$ is even, then $\tilde{E}^n_k(Z,s)$ is holomorphic for $k+2\Re(s)>[n/2]$ except for a possible simple pole at $k+2\Re(s)=n/2+1$.

\noindent{\rm (iii)} If\/ $k<n/2$ then\/ $\tilde{E}^n_k(Z,s)$ is holomorphic for\/ $k+2\Re(s)>[n/2]$ except for possible simple poles at
\[ \left[\frac{n+3}{2}\right]\leq k+2s\leq n-k+1, \quad k+2s\in\mathbb{Z}. \]
}

\bigskip
\noindent{\bf 3. Links and differential operators}

\bigskip
\noindent{\bf 3.1. Links, homogeneous polynomials and differential operators}

 Fix an index set $I$. A non-ordered pair $(\alpha_1,\alpha_2)$ is called a {\sl link\/} if $\alpha_1$, $\alpha_2\in I$ are distinct. Let the set $\mathcal{L}(I)$ of links be
\[ \mathcal{L}(I) := \{ \{(\alpha_1,\alpha_2),(\alpha_3,\alpha_4)\ldots,(\alpha_{2r-1},\alpha_{2r})\} \mid \alpha_1,\alpha_2,\ldots,\alpha_{2r}\in I, \quad \alpha_i\neq\alpha_j (i\neq j) \quad\text{for some $r$} \}. \]
For $L=\{(\alpha_1,\alpha_2),(\alpha_3,\alpha_4),\ldots,(\alpha_{2r-1},\alpha_{2r})\}\in\mathcal{L}(I)$, $\overline{L}$ denotes the set $\{\alpha_1,\alpha_2,\ldots,\alpha_{2r}\}$. We remark $\#L=r$ and $\#\overline{L}=2r$.

 For any link $(\alpha_1,\alpha_2)$, let $X^{\alpha\beta}$ be an indeterminate. Define $\mathbb{C}[X;I]$ to be the polynomial ring generated by
\[ \{X^{\alpha_1\alpha_2} \mid \text{$(\alpha_1,\alpha_2)$ is a link} \}. \]
Furthermore for indeterminates $Y^{\alpha_1\alpha_2}$, define $\mathbb{C}[X,Y;I]$ to be the polynomial ring generated by
\[ \{X^{\alpha_1\alpha_2}, Y^{\alpha_1\alpha_2} \mid \text{$(\alpha_1,\alpha_2)$ is a link} \}. \]
For $L=\{(\alpha_1,\alpha_2),(\alpha_3,\alpha_4),\ldots,(\alpha_{2r-1},\alpha_{2r})\}\in\mathcal{L}(I)$, let
\[ X^L := X^{\alpha_1\alpha_2} X^{\alpha_3\alpha_4}\ldots X^{\alpha_{2r-1}\alpha_{2r}}\in\mathbb{C}[X;I]. \]
Let
\[ \mathbb{C}[X;\mathcal{L}(I)] := \langle X^L \mid L\in\mathcal{L}(I)\rangle_{\mathbb{C}},\quad \text{$\langle\quad\rangle_{\mathbb{C}}$ denotes the $\mathbb{C}$-linear span).} \]
$P(X)\in\mathbb{C}[X;\mathcal{L}(I)]$ is called a {\sl homogeneous polynomial} if there exists a subset $J$ of $I$ such that
\[ P(X) = \sum_{\overline{L}=J} c(L) X^L, \quad c(L)\in\mathbb{C}. \]
This subset $J$ is called the {\sl underlying set}. $\#J/2$ is the {\sl degree of\/ $P$} and denoted by $\deg P$.

 Let $\delta$, $\varepsilon$, $\Delta$ and $\Epsilon$ be indeterminates. For any link $(\alpha_1,\alpha_2)$, define a linear differential operator $\partial^{\alpha_1\alpha_2}$ satisfying:

(D1) $\partial^{\alpha_1\alpha_2}(fg) = (\partial^{\alpha_1\alpha_2}f)g + f(\partial^{\alpha_1\alpha_2}g)$,

(D2) $\partial^{\alpha_1\alpha_2}(\delta^{-k}\left|\delta\right|^{-2s}\varepsilon^s) = \delta^{-k}\left|\delta\right|^{-2s}\varepsilon^s \bigl((-k-s)\Delta^{\alpha_1\alpha_2}+s\Epsilon^{\alpha_1\alpha_2}\bigr)$,

(D3) $\partial^{\alpha_1\alpha_2} \Delta^{\alpha_3\alpha_4} = -\frac{1}{2} (\Delta^{\alpha_1\alpha_3}\Delta^{\alpha_2\alpha_4} + \Delta^{\alpha_1\alpha_4}\Delta^{\alpha_2\alpha_3})$ if $\alpha_1$, $\ldots$, $\alpha_4$ are distinct,

(D4) $\partial^{\alpha_1\alpha_2} \Epsilon^{\alpha_3\alpha_4} = -\frac{1}{2} (\Epsilon^{\alpha_1\alpha_3}\Epsilon^{\alpha_2\alpha_4} + \Epsilon^{\alpha_1\alpha_4}\Epsilon^{\alpha_2\alpha_3})$ if $\alpha_1$, $\ldots$, $\alpha_4$ are distinct.

\noindent Furthermore, for $L=\{(\alpha_1,\alpha_2),(\alpha_3,\alpha_4),\ldots,(\alpha_{2r-1},\alpha_{2r})\}\in\mathcal{L}(I)$, define
\[ \partial^L := \partial^{\alpha_1\alpha_2}\partial^{\alpha_3\alpha_4}\ldots\partial^{\alpha_{2r-1}\alpha_{2r}}. \]

 Then for $L_0\in\mathcal{L}(I)$, we have the following form:
\[ \partial^{L_0}(\delta^{-k}\left|\delta\right|^{-2s}\varepsilon^s) = \delta^{-k}\left|\delta\right|^{-2s}\varepsilon^s \sum c(L_1,L_2) (\Delta-\Epsilon)^{L_1} \Epsilon^{L_2}, \tag{3.1} \]
where the summation runs over
\[ \{(L_1,L_2) \mid L=L_1\sqcup L_2\in\mathcal{L},\quad\overline{L}=I,\quad L_1,L_2\in\mathcal{L}\sqcup\{\emptyset\}\}. \]

\bigskip
\noindent{\bf 4. Differential operators for Siegel modular forms} \par

 Let $W_j$ be a $j$-dimensional vector space
\[ W_j := \mathbb{C}e_1\oplus\mathbb{C}e_2\oplus\cdots\oplus\mathbb{C}e_j, \]
where $e_1$, $e_2$, $\ldots$, $e_j$ are indeterminates. Let $T^l(W_j)$ be the $l$th tensor product of $W_j$ for a positive integer $l$, i.e.,
\[ T^l(W_j):=\underbrace{W_j\otimes\cdots\otimes W_j}_{\text{$l$-times}}, \]
and $\rho_j'$ the standard representation of $GL(j,\mathbb{C})$ on $T^l(W_j)$. Fix a positive integer $m$ and the Young symmetrizer $c$ for an $m$-tuple $(\lambda_1,\lambda_2,\ldots,\lambda_m)$. Here $(\lambda_1,\lambda_2,\ldots,\lambda_m)\in\mathbb{Z}^m$, $\lambda_1\geq\lambda_2\geq\ldots\geq\lambda_m>0$ and $\lambda_1+\lambda_2+\cdots+\lambda_m=l$ (see Weyl [21, Chapter 4]). Put $\rho_j:=\det^k\otimes \rho_j'$ for an integer $k$ and $V_j:=c(T^l(W_j))$. Then $(\rho_j,V_j)$ is an irreducible representation. Note that $V_j=\{0\}$ if and only if $j<m$.

 Let ${W_j}^*$ and ${W_j}_*$ be copies of $W_j$, i.e.,
\[ {W_j}^* := \mathbb{C}{e_1}^*\oplus\mathbb{C}{e_2}^*\oplus\cdots\oplus\mathbb{C}{e_j}^*,\quad {W_j}_* := \mathbb{C}{e_1}_*\oplus\mathbb{C}{e_2}_*\oplus\cdots\oplus\mathbb{C}{e_j}_*, \]
where ${e_1}^*$, ${e_2}^*$, $\ldots$, ${e_j}^*$, ${e_1}_*$, ${e_2}_*$, $\ldots$, ${e_j}_*$ are indeterminates. For $w\in W_j$ and $v\in T^l(W_j)$, $w^*\in{W_j}^*$, $w_*\in{W_j}_*$, $v^*\in T^l({W_j}^*)$ and $v_*\in T^l({W_j}_*)$ are defined in the obvious way. For ${W_j}^*$ and the same $l$, $m$, $(\lambda_1,\lambda_2,\ldots,\lambda_m)$ and $k$, define $c^*$, ${\rho_j}^*$ and ${V_j}^*$ as the same as above. Similarly, for ${W_j}_*$, define $c_*$, ${\rho_j}_*$ and ${V_j}_*$.

 Let $e_j^{(\alpha)}$ be an indeterminate for $\alpha\in I=I^*\cup I_*$ and a positive integer $j$. For a symmetric matrix $A$ of size $p+q$ and positive integers $a$, $b$, define
\begin{align*}
    A^{ab} &:= (e_1^{(a^*)},\ldots,e_p^{(a^*)},0,\ldots,0) \,A\, {}^t\!(e_1^{(b^*)},\ldots,e_p^{(b^*)},0,\ldots,0),
\\  A^a_b  &:= (e_1^{(a^*)},\ldots,e_p^{(a^*)},0,\ldots,0) \,A\, {}^t\!(0,\ldots,0,e_1^{(b_*)},\ldots,e_q^{(b_*)}),
\\  A_{ab} &:= (0,\ldots,0,e_1^{(a_*)},\ldots,e_q^{(a_*)}) \,A\, {}^t\!(0,\ldots,0,e_1^{(b_*)},\ldots,e_q^{(b_*)}).
\end{align*}
We identify $e_{i_1}^{(1^*)}e_{i_2}^{(2^*)}\ldots e_{i_l}^{(l^*)}e_{j_1}^{(1_*)}e_{j_2}^{(2_*)}\ldots e_{j_l}^{(l_*)}$ with $e_{i_1}^*\otimes e_{i_2}^*\otimes\cdots\otimes e_{i_l}^*\otimes{e_{j_1}}_*\otimes{e_{j_2}}_*\otimes\cdots\otimes{e_{j_l}}_*$.

 For $g=\begin{pmatrix} \mathfrak{A} & \mathfrak{B} \\ \mathfrak{C} & \mathfrak{D} \end{pmatrix}\in Sp(p+q,\mathbb{R})$ and $\mathfrak{Z}=(z_{\mu\nu})_{1\leq\mu,\nu\leq p+q}\in\mathfrak{H}_{p+q}$, put
\[ \partial := \left( \frac{1+\delta_{\mu\nu}}{2} \frac{\partial}{\partial z_{\mu\nu}} \right)_{1\leq\mu,\nu\leq p+q}, \quad \delta := \delta(g,\mathfrak{Z}) := \det(\mathfrak{C}\mathfrak{Z}+\mathfrak{D}), \quad \Delta := \Delta(g,\mathfrak{Z}) := (\mathfrak{C}\mathfrak{Z}+\mathfrak{D})^{-1}\mathfrak{C}, \]
\[ \varepsilon := \det(\Im(\mathfrak{Z})) \quad\text{and}\quad \Epsilon := \frac{1}{2i}(\Im(\mathfrak{Z}))^{-1}, \]
where $\delta_{\mu\nu}$ is Kronecker's delta. Then $\partial$, $\delta$, $\varepsilon$, $\Delta$ and $\Epsilon$ satisfy (D1)--(D4) in Section 3.

 Let $\mathfrak{P}(X)$ be a homogeneous polynomial with the underlying set $I$ and pluri-harmonic for $\Delta^*(d)$ and $\Delta_*(d)$. Put $X:=\begin{pmatrix} X_1 \\ X_2 \end{pmatrix} {\Bigm.}^t\!\begin{pmatrix} X_1 \\ X_2 \end{pmatrix}$, where $X_1\in M(p,d,\mathbb{C})$ and $X_2\in M(q,d,\mathbb{C})$. Then $\mathfrak{P}(X)$ is pluri-harmonic for each $X_1$ and $X_2$, i.e., for $X_1=(x_{\mu\nu}^{(1)})$ and $X_2=(x_{\mu\nu}^{(2)})$,
\[ \sum_{\kappa=1}^d\frac{\partial}{\partial x_{\mu\kappa}^{(1)}}\frac{\partial}{\partial x_{\nu\kappa}^{(1)}}\mathfrak{P}(X) = 0 \quad\text{and}\quad \sum_{\kappa=1}^d\frac{\partial}{\partial x_{\mu\kappa}^{(2)}}\frac{\partial}{\partial x_{\nu\kappa}^{(2)}}\mathfrak{P}(X) = 0 \quad\text{for any $\mu$, $\nu$}. \]
Using the same method as in [7], we have:

\bigskip
{\bf Theorem 4.1.} (cf. Ibukiyama [7, Theorem 1])
{\sl
Let\/ $\mathfrak{P}(X)$ be a homogeneous polynomial with the underlying set\/ $I$ and pluri-harmonic for\/ $\Delta^*(2k)$ and\/ $\Delta_*(2k)$. For a\/ $C^\infty$-function $f\colon\mathfrak{H}_{p+q}\to\mathbb{C}$, $g_1\in Sp(p,\mathbb{R})$ and\/ $g_2\in Sp(q,\mathbb{R})$, we have
\[ ((\mathfrak{P}(\partial)f)|_{\rho_p^*} g_1|_{\rho_{q*}} g_2)|_{\mathfrak{Z}=\mathfrak{Z}_0} = (\mathfrak{P}(\partial)(f|_k g_1^* g_{2*}))|_{\mathfrak{Z}=\mathfrak{Z}_0}, \]
where $g_1^*:=\begin{pmatrix} A_1 & 0 & B_1 & 0 \\ 0 & 1_q & 0 & 0 \\ C_1 & 0 & D_1 & 0 \\ 0 & 0 & 0 & 1_q \end{pmatrix}$ and\/ ${g_2}_*:=\begin{pmatrix} 1_p & 0 & 0 & 0 \\ 0 & A_2 & 0 & B_2 \\ 0 & 0 & 1_p & 0 \\ 0 & C_2 & 0 & D_2 \end{pmatrix}$ for\/ $g_1=\begin{pmatrix} A_1^{(p)} & B_1^{(p)} \\ C_1^{(p)} & D_1^{(p)} \end{pmatrix}$ and\/ $g_2=\begin{pmatrix} A_2^{(q)} & B_2^{(q)} \\ C_2^{(q)} & D_2^{(q)} \end{pmatrix}$, and\/ $\mathfrak{Z}_0:=\begin{pmatrix} Z^{(p)} & 0 \\ 0 & W^{(q)} \end{pmatrix}$.
}

\bigskip

Using the above theorem, we prove the following proposition.

\bigskip
{\bf Proposition 4.2.}
{\sl
If\/ $\mathfrak{P}(X)$ is a homogeneous polynomial in\/ $X$ with the underlying set $I$ and pluri-harmonic for $\Delta^*(2k)$ and $\Delta_*(2k)$, then there exists a homogeneous polynomial\/ $\mathfrak{Q}(X,s)$ of\/ $X$ such that
\[ \bigl(\mathfrak{P}(\partial) (\delta^{-k}\left|\delta\right|^{-2s}\varepsilon^s)\bigr)|_{\mathfrak{Z}=\mathfrak{Z}_0} = \bigl(\delta^{-k}\left|\delta\right|^{-2s}\varepsilon^s \mathfrak{Q}(\Delta-\Epsilon,s)\bigr)|_{\mathfrak{Z}=\mathfrak{Z}_0}. \]
}
\bigskip

Proof.

Without loss of generality, $n:=p=q$ and $n\geq l$. By Theorem 4.1, for $g\in Sp(2n,\mathbb{R})$, $g_1$, $g_2\in Sp(n,\mathbb{R})$,
\[ (\mathfrak{P}(\partial)(\det(\Im(\mathfrak{Z}))^s|_kg)|_{\rho_n^*}g_1|_{\rho_{n*}}g_2)|_{\mathfrak{Z}=\mathfrak{Z}_0} = (\mathfrak{P}(\partial)(\det(\Im(\mathfrak{Z}))^s|_k gg_1^*g_{2*})|_{\mathfrak{Z}=\mathfrak{Z}_0}. \tag{4.1} \]
Since (3.1), we can write $\mathfrak{P}(\partial)(\delta^{-k}\left|\delta\right|^{-2s}\varepsilon^s)$ in the form
\[ \mathfrak{P}(\partial) (\delta^{-k}\left|\delta\right|^{-2s}\varepsilon^s) = \delta^{-k}\left|\delta\right|^{-2s}\varepsilon^s \sum_{(L_1,L_2)} c(L_1,L_2) (\Delta-\Epsilon)^{L_1}\Epsilon^{L_2}, \tag{4.2} \]
where the summation runs over
\[ \{(L_1,L_2) \mid L=L_1\sqcup L_2\in\mathcal{L},\quad\overline{L}=I,\quad L_1,L_2\in\mathcal{L}\sqcup\{\emptyset\}\}. \]
Substituting (4.2) into (4.1), the left-hand side of (4.1) is equal to
\begin{align*}
    & \delta(gg_1^*g_{2*},\mathfrak{Z}_0)^{-k}\left|\delta(gg_1^*g_{2*},\mathfrak{Z}_0)\right|^{-2s}\det(\Im(\mathfrak{Z}_0))^s
\\  & \quad \cdot \sum_{(L_1,L_2)} c(L_1,L_2) \left(-\frac{1}{2i}j(gg_1^*g_{2*},\mathfrak{Z}_0)^{-1}(\Im(gg_1^*g_{2*}\langle\mathfrak{Z}_0\rangle))^{-1}{}^t\!j(gg_1^*g_{2*},\mathfrak{Z}_0)^{-1}\right)^{L_1}
\\  & \quad\qquad \cdot \left(\frac{1}{2i}j(g_1^*g_{2*},\mathfrak{Z}_0)^{-1}(\Im(g_1^*g_{2*}\langle\mathfrak{Z}_0\rangle))^{-1}{}^t\!j(g_1^*g_{2*},\mathfrak{Z}_0)^{-1}\right)^{L_2}
\end{align*}
and the right-hand side is equal to
\begin{align*}
    & \delta(gg_1^*g_{2*},\mathfrak{Z}_0)^{-k}\left|\delta(gg_1^*g_{2*},\mathfrak{Z}_0)\right|^{-2s}\det(\Im(\mathfrak{Z}_0))^s
\\  & \quad \cdot \sum_{(L_1,L_2)} c(L_1,L_2) \left(-\frac{1}{2i}j(gg_1^*g_{2*},\mathfrak{Z}_0)^{-1}(\Im(gg_1^*g_{2*}\langle\mathfrak{Z}_0\rangle))^{-1}{}^t\!j(gg_1^*g_{2*},\mathfrak{Z}_0)^{-1}\right)^{L_1}\left(\frac{1}{2i}(\Im(\mathfrak{Z}_0))^{-1}\right)^{L_2}.
\end{align*}
Therefore we have
\begin{align*}
    \sum_{(L_1,L_2)} & c(L_1,L_2)\left(-\frac{1}{2i}j(gg_1^*g_{2*},\mathfrak{Z}_0)^{-1}(\Im(gg_1^*g_{2*}\langle\mathfrak{Z}_0\rangle))^{-1}{}^t\!j(gg_1^*g_{2*},\mathfrak{Z}_0)^{-1}\right)^{L_1} \tag{4.3}
\\                   & \cdot \left\{ \left(\frac{1}{2i}j(g_1^*g_{2*},\mathfrak{Z}_0)^{-1}(\Im(g_1^*g_{2*}\langle\mathfrak{Z}_0\rangle))^{-1}{}^t\!j(g_1^*g_{2*},\mathfrak{Z}_0)^{-1}\right)^{L_2} - \left(\frac{1}{2i}(\Im(\mathfrak{Z}_0))^{-1}\right)^{L_2} \right\} = 0.
\end{align*}
Clearly the terms of $L_2=\emptyset$ vanish. And if $L_2\neq\emptyset$ and $(\alpha,\beta)\in\mathcal{L}^*_*$ for some $(\alpha,\beta)\in L_2$, then
\[ \left(\frac{1}{2i}j(g_1^*g_{2*},\mathfrak{Z}_0)^{-1}(\Im(g_1^*g_{2*}\langle\mathfrak{Z}_0\rangle))^{-1}{}^t\!j(g_1^*g_{2*},\mathfrak{Z}_0)^{-1}\right)^{L_2} - \left(\frac{1}{2i}(\Im(\mathfrak{Z}_0))^{-1}\right)^{L_2} \]
vanishes. Therefore the summation of (4.3) runs over
\[ \mathcal{L}':=\{(L_1,L_2)\mid\text{$L_2\neq\emptyset$ and $(\alpha,\beta)\notin\mathcal{L}^*_*$ for all $(\alpha,\beta)\in L_2$}\}. \]
It suffices to prove that $c(L_1,L_2)=0$ for $(L_1,L_2)\in\mathcal{L}'$.

Put $\mathfrak{Z}_0=i1_{2n}$ and $g$, $g_1^*g_{2*}\in K$, where
\[ K:=\{g\in Sp(2n,\mathbb{R}) \mid g\langle i1_{2n}\rangle = i1_{2n}\}. \]
Then the left-hand side of (4.3) is equal to
\begin{align*}
    \sum_{(L_1,L_2)} & c(L_1,L_2) \left(-\frac{1}{2i}j(gg_1^*g_{2*},i1_{2n})^{-1} {}^t\!j(gg_1^*g_{2*},i1_{2n})^{-1}\right)^{L_1} \tag{4.4}
\\                   & \quad \cdot \left\{ \left(\frac{1}{2i}j(g_1^*g_{2*},i1_{2n})^{-1}{}^t\!j(g_1^*g_{2*},i1_{2n})\right)^{L_2} - \left(\frac{1}{2i}1_{2n}\right)^{L_2} \right\}.
\end{align*}
Let
\[ g = \begin{pmatrix} \Re(u) & -\Im(u) \\ \Im(u) & \Re(u) \end{pmatrix} \quad\text{with}\quad u = \begin{pmatrix} a_1 & & & b_1 & & \\ & \ddots & & & \ddots & \\ & & a_n & & & b_n \\ -\overline{b_1} & & & \overline{a_1} & & \\ & \ddots & & & \ddots & \\ & & -\overline{b_n} & & & \overline{a_n} \end{pmatrix}^{-1}\in U(2n), \]
and
\[ g_1^*g_{2*} = \begin{pmatrix} \Re(v) & -\Im(v) \\ \Im(v) & \Re(v) \end{pmatrix} \quad\text{with}\quad v = \begin{pmatrix} e^{i\theta_1} & & \\ & \ddots & \\ & & e^{i\theta_{2n}} \end{pmatrix}^{-1} \]
for $\theta_1$, $\ldots$, $\theta_{2n}\in\mathbb{R}$. Here $U(2n)$ denotes the unitary group of size $2n$. Then (4.4) is equal to
\[ \sum_{(L_1,L_2)} c(L_1,L_2) \left(-\frac{1}{2i}v^{-1}u^{-1}{}^t\!u^{-1}{}^t\!v^{-1}\right)^{L_1} \left\{ \left(\frac{1}{2i}v^{-1}{}^t\!v^{-1}\right)^{L_2} - \left(\frac{1}{2i}1_{2n}\right)^{L_2} \right\}. \]
Since $n\geq l$, for $L=\{(\alpha_1,\alpha_2),(\alpha_3,\alpha_4),\ldots,(\alpha_{2l-1},\alpha_{2l})\}\in\mathcal{L}$, the coefficient of $e_1^{(\alpha_1)}e_1^{(\alpha_2)}e_2^{(\alpha_3)}e_2^{(\alpha_4)}\ldots$ $e_l^{(\alpha_{2l-1})}e_l^{(\alpha_{2l})}$ is
\[ \sum c(L_1,L_2) \left(\prod_{(\alpha_{2j-1},\alpha_{2j})\in L_1} c_j \right) \left\{ \left(\prod_{(\alpha_{2j-1},\alpha_{2j})\in L_2} d_j \right) - \left(\frac{1}{2i}\right)^{\# L_2} \right\}. \]
Here the summation runs over
\[ \mathcal{L}'':=\{(L_1,L_2)\in\mathcal{L}'\mid L_1\sqcup L_2 = \{(\alpha_1,\alpha_2),\ldots,(\alpha_{2l-1},\alpha_{2l})\}, \]
and we put
\[ c_j:=
    \begin{cases}
     -\frac{1}{2i} e^{2i\theta_j}(a_j^2+b_j^2), &\quad\text{if $(\alpha_{2j-1},\alpha_{2j})\in\mathcal{L}^{**}$}
\\   -\frac{1}{2i} e^{i(\theta_j+\theta_{n+j})}(\overline{a_j}b_j-a_j\overline{b_j}), &\quad\text{if $(\alpha_{2j-1},\alpha_{2j})\in\mathcal{L}^*_*$}
\\   -\frac{1}{2i} e^{2i\theta_{n+j}}(\overline{a_j}^2+\overline{b_j}^2), &\quad\text{if $(\alpha_{2j-1},\alpha_{2j})\in\mathcal{L}_{**}$}
    \end{cases} \]
and
\[ d_j:=
    \begin{cases}
     \frac{1}{2i} e^{2i\theta_j},     &\quad\text{if $(\alpha_{2j-1},\alpha_{2j})\in\mathcal{L}^{**}$}
\\   \frac{1}{2i} e^{2i\theta_{n+j}}, &\quad\text{if $(\alpha_{2j-1},\alpha_{2j})\in\mathcal{L}_{**}$}
    \end{cases}. \]
Since $n\geq l$,
\[ \left\{ \left(\prod_{(\alpha_{2j-1},\alpha_{2j})\in L_1} c_j \right) \left\{ \left(\prod_{(\alpha_{2j-1},\alpha_{2j})\in L_2} d_j \right) - \left(\frac{1}{2i}\right)^{\# L_2} \right\} \Biggm| (L_1,L_2)\in\mathcal{L}'' \right\} \]
is linearly independent. Therefore we obtain $c(L_1,L_2)=0$ for $(L_1,L_2)\in\mathcal{L}''$. Thus Proposition 4.2 is proved. \qed

\bigskip
\noindent{\bf 5. Pullback formula} \par

\bigskip
{\bf Proposition 5.1.}
{\sl
 Let\/ $k$ be even and\/ $k+2\Re(s)>p+q+1$. If\/ $\mathfrak{P}(X)$ is a homogeneous polynomial in\/ $X$ with the underlying set\/ $I$ and pluri-harmonic for\/ $\Delta^*(2k)$ and\/ $\Delta_*(2k)$, then
\[ (\mathfrak{P}(\partial)E^{p+q}_k)(\begin{pmatrix} Z^{(p)} & 0 \\ 0 & W^{(q)} \end{pmatrix},s) = \sum_{r=0}^{\min(p,q)} \sum_{T\in\mathbb{T}^{(r)}} \mathcal{P}_r(Z,W,T,s), \]
where
\begin{align*}
    \mathcal{P}_r(Z,W&,T,s) := \sum_{g_2''\in\Gamma_{q,r}\backslash\Gamma_q} \sum_{g_1''\in\Gamma_{p,r}\backslash\Gamma_p} \sum_{g_2'\in\Gamma_r(T)\backslash\Gamma_r} \sum_{g_1'\in\Gamma_r}
\\  & \cdot \bigl\{\det(\Im(Z))^s \det(\Im(W))^s \left|\det(1_p-\tilde{T}W{}^t\!\tilde{T}Z)\right|^{-2s}
\\  & \cdot \det(1_p-\tilde{T}W{}^t\!\tilde{T}Z)^{-k} \mathfrak{Q}(P,s) \bigr\} |_{\rho_p^*}\tilde{g_1'}|_{\rho_{q*}}\tilde{g_2'}|_{\rho_p^*}g_1''|_{\rho_{q*}}g_2'',
\end{align*}
\[ \Gamma_r(T):=\begin{pmatrix} T & 0 \\ 0 & T^{-1}\end{pmatrix}^{-1}\Gamma_r\begin{pmatrix} T & 0 \\ 0 & T^{-1}\end{pmatrix} \cap \Gamma_r, \]
\[ \tilde{T} := \begin{pmatrix} T^{(r)} & 0 \\ 0 & 0 \end{pmatrix}, \quad\text{and}\quad g_{\tilde{T}} := \begin{pmatrix} 1_p & 0 & 0 & 0 \\ 0 & 1_q & 0 & 0 \\ 0 & \tilde{T} & 1_p & 0 \\ {}^t\!\tilde{T} & 0 & 0 & 1_q \end{pmatrix}, \]
\[ P := -\frac{1}{2i} \begin{pmatrix} 1_p-\tilde{T}W{}^t\!\tilde{T}Z & 0 \\ 0 & 1_q-{}^t\!\tilde{T}Z\tilde{T}W \end{pmatrix}^{-1} \begin{pmatrix} (1_p-\tilde{T}W{}^t\!\tilde{T}\overline{Z})(\Im(Z))^{-1} & -2i\tilde{T} \\ -2i{}^t\!\tilde{T} & (1_q-{}^t\!\tilde{T}Z\tilde{T}\overline{W})(\Im(W))^{-1} \end{pmatrix}. \]
}
\bigskip

 Proof. By Garrett [6] (see also Mizumoto [15, Appendix B]), the left coset $P_{p+q,0}\backslash\Gamma_{p+q}$ has a complete system of representatives $g_{\tilde{T}} \tilde{g_1'}^* \tilde{g_2'}_* {g_1''}^* {g_2''}_*$ with
\[ g_{\tilde{T}} = \begin{pmatrix} 1_p & 0 & 0 & 0 \\ 0 & 1_q & 0 & 0 \\ 0 & \tilde{T} & 1_p & 0 \\ {}^t\!\tilde{T} & 0 & 0 & 1_q \end{pmatrix} \quad \text{for $T\in\mathbb{T}^{(r)}$ \quad ($0\leq r\leq\min(p,q)$)}, \]
\[ g_1'\in\Gamma_r,\quad g_1''\in\Gamma_{p,r}\backslash\Gamma_p,\quad g_2'\in\Gamma_r(T)\backslash\Gamma_r,\quad g_2''\in\Gamma_{q,r}\backslash\Gamma_q. \]
We put $g_1:=\tilde{g_1'}g_1''$, $g_2:=\tilde{g_2'}g_2''$ and $\begin{pmatrix} * & * \\ \mathfrak{C}^{(p+q)} & \mathfrak{D}^{(p+q)} \end{pmatrix}:=g_{\tilde{T}}g_1^*g_{2*}$. Then, by Proposition 4.2, we have
\begin{align*}
    & \mathfrak{P}(\partial)(\det(\Im(\mathfrak{Z}))^s \det(\mathfrak{C}\mathfrak{Z}+\mathfrak{D})^{-k} \left|\det(\mathfrak{C}\mathfrak{Z}+\mathfrak{D})\right|^{-2s}) |_{\mathfrak{Z}=\mathfrak{Z}_0} \tag{5.1}
\\  & = \det(\Im(\mathfrak{Z}_0))^s \det(\mathfrak{C}\mathfrak{Z}_0+\mathfrak{D})^{-k} \left|\det(\mathfrak{C}\mathfrak{Z}_0+\mathfrak{D})\right|^{-2s} \cdot \mathfrak{Q}((\mathfrak{C}\mathfrak{Z}_0+\mathfrak{D})^{-1}\mathfrak{C}-\frac{1}{2i}(\Im(\mathfrak{Z}_0))^{-1},s).
\end{align*}
Since
\[ \mathfrak{C}\mathfrak{Z}_0+\mathfrak{D}=\begin{pmatrix} 1_p & \tilde{T}g_2\langle W\rangle \\ {}^t\!\tilde{T}g_1\langle Z\rangle & 1_q\end{pmatrix} \begin{pmatrix} j(g_1,Z) & 0 \\ 0 & j(g_2,W) \end{pmatrix} \]
and
\[ (\mathfrak{C}\mathfrak{Z}_0+\mathfrak{D})^{-1}\mathfrak{C}-\frac{1}{2i}(\Im(\mathfrak{Z}_0))^{-1} = -\frac{1}{2i}(\mathfrak{C}\mathfrak{Z}_0+\mathfrak{D})^{-1}(\mathfrak{C}\overline{\mathfrak{Z}_0}+\mathfrak{D})(\Im(\mathfrak{Z}_0))^{-1}, \]
we obtain
\[ \det(\mathfrak{C}\mathfrak{Z}_0+\mathfrak{D}) = \det(j(g_1,Z))\det(j(g_2,W))\det(1_p-\tilde{T}g_2\langle W\rangle{}^t\!\tilde{T}g_1\langle Z\rangle) \]
and
\begin{align*}
    & (\mathfrak{C}\mathfrak{Z}_0+\mathfrak{D})^{-1}\mathfrak{C}-\frac{1}{2i}(\Im(\mathfrak{Z}_0))^{-1}
\\  & = -\frac{1}{2i}\begin{pmatrix} j(g_1,Z) & 0 \\ 0 & j(g_2,W) \end{pmatrix}^{-1} \begin{pmatrix} 1_p-\tilde{T}g_2\langle W\rangle{}^t\!\tilde{T}g_1\langle Z\rangle & 0 \\ 0 & 1_q-{}^t\!\tilde{T}g_1\langle Z\rangle\tilde{T}g_2\langle W\rangle \end{pmatrix}^{-1}
\\  & \qquad\cdot \begin{pmatrix} (1_p-\tilde{T}g_2\langle W\rangle{}^t\!\tilde{T}g_1\langle\overline{Z}\rangle)(\Im(g_1\langle Z\rangle))^{-1} & -2i\tilde{T} \\ -2i{}^t\!\tilde{T} & (1_q-{}^t\!\tilde{T}g_1\langle Z\rangle\tilde{T}g_2\langle\overline{W}\rangle)(\Im(g_2\langle W\rangle))^{-1} \end{pmatrix}
\\  & \qquad\cdot {\Bigm.}^t\!\begin{pmatrix} j(g_1,Z) & 0 \\ 0 & j(g_2,W) \end{pmatrix}^{-1}.
\end{align*}
Therefore (5.1) is equal to
\[ \det(\Im(Z))^s \det(\Im(W))^s \left|\det(1_p-\tilde{T}W{}^t\!\tilde{T}Z)\right|^{-2s}\det(1_p-\tilde{T}W{}^t\!\tilde{T}Z)^{-k}\mathfrak{Q}(P,s)|_{\rho_p^*}g_1|_{\rho_{q*}}g_2. \]
Thus Proposition 5.1 is proved. \qed

\bigskip
{\bf Theorem 5.2.}
{\sl
 Let\/ $k$ be even and\/ $k+2\Re(s)>p+q+1$. If\/ $V_p^*\otimes V_{q*}$-valued polynomial $\mathfrak{P}(X)$ is homogeneous in\/ $X$ with underlying set\/ $I$ and pluri-harmonic for\/ $\Delta^*(2k)$ and\/ $\Delta_*(2k)$, then for any eigenform\/ $f\in\mathfrak{S}_{\rho_q}$,
\begin{align*}
    & \left( f_*,(\mathfrak{P}(\partial)E^{p+q}_k)\Bigl(\begin{pmatrix} -\overline{Z}^{(p)} & 0 \\ 0 & * \end{pmatrix},\overline{s} \Bigr) \right)
\\  & \qquad\qquad = 2^{q(q+1-2s)+1}\, (2^{-1}i)^{\left|\rho_q\right|}\, c(s,\rho_q)\, D(k+2s,f)\, [f]^{\rho_p}_{\rho_q}(Z,s)^*.
\end{align*}
Here a\/ $\mathbb{C}$-valued function\/ $c(s,\rho_q)$ satisfies as follows:
\[ \int_{S_q}\langle\rho_{q*}(1_q-\overline{S}S)v_*,\mathfrak{Q}(R,\overline{s})\rangle \det(1_q-\overline{S}S)^{s-q-1}\, dS = c(s,\rho_q) v^* \quad\text{for any $v\in V_q$}, \]
with\/ $R:=\displaystyle{-\frac{1}{2i}}\begin{pmatrix} S & -2i 1_q \\ -2i 1_q & 2^2\overline{S}(1_q-S\overline{S})^{-1} \end{pmatrix}$.
}
\bigskip

 Proof of Theorem 5.2. First we consider that the Petersson inner product $(f_*,\mathcal{P}_r(-\overline{Z},*,T,\overline{s}))$. By the same reason as that Klingen [10, Satz 2],
\[ (f_*,\mathcal{P}_r(-\overline{Z},*,T,\overline{s}))=0 \quad\text{for $r<q$}. \tag{5.2} \]
Therefore we only consider that $(f_*,\mathcal{P}_q(-\overline{Z},*,T,\overline{s}))$.

 Now, we have
\begin{align*}
    \mathcal{P}_q(Z,&W,T,s) = \sum_{g_1''\in\Gamma_{p,q}\backslash\Gamma_p} \sum_{g_2'\in\Gamma_q(T)\backslash\Gamma_q} \bigl\{ \det(T)^{-k-2s} \left(\frac{\det(\Im(Z))}{\det(\Im(Z_1))}\right)^s
\\             & \cdot \mathcal{P}(Z,W,s) \bigr\} \Bigm|_{\rho_{q*}} \begin{pmatrix} T & 0 \\ 0 & T^{-1} \end{pmatrix} g_2' \Bigm|_{\rho_p^*} g_1'',
\end{align*}
where
\begin{align*}
    \mathcal{P}(Z,&W,s) = \sum_{g_1'\in\Gamma_q} \bigl\{ \det(\Im(Z_1))^s \det(\Im(W))^s \left|\det(Z_1+W)\right|^{-2s}
\\  & \cdot \det(Z_1+W)^{-k} \mathfrak{Q}(Q(Z,W),s) \bigr\}|_{\rho_p^*}\tilde{g}_1'
\end{align*}
with
\begin{align*}
    & Q(Z,W)
\\  & = -\frac{1}{2i} \begin{pmatrix} \begin{pmatrix} Z_1+W & Z_2 \\ 0 & 1_{p-q} \end{pmatrix}^{-1} \begin{pmatrix} \overline{Z_1}+W & \overline{Z_2} \\ 0 & 1_{p-q} \end{pmatrix} (\Im(Z))^{-1} & -2i \begin{pmatrix} (Z_1+W)^{-1} \\ 0 \end{pmatrix} \\ -2i \begin{pmatrix} (Z_1+W)^{-1} & 0 \end{pmatrix} & (Z_1+W)^{-1}(Z_1+\overline{W})(\Im(W))^{-1} \end{pmatrix}.
\end{align*}
Since the Hecke operators are Hermitian operators and $f$ is an eigenform, we have
\[ (f_*,\mathcal{P}_q(-\overline{Z},*,T,\overline{s})) = \lambda(T,f) \det(T)^{-k-2s} \sum_{g_1''\in\Gamma_{p,q}\backslash\Gamma_p} \left(\frac{\det(\Im(Z))}{\det(\Im(Z_1))}\right)^s (f_*,\mathcal{P}(-\overline{Z},*,\overline{s}))\Bigm|_{\rho_p^*} g_1''. \tag{5.3} \]

{\bf Lemma 5.3.}
{\sl
\[ (f_*,\mathcal{P}(-\overline{Z},*,\overline{s})) = 2^{q(q+1-2s)+1}\, (2^{-1}i)^{\left|\rho_q\right|} \cdot c(s,\rho_q) \cdot f(Z)^*. \]
}

 Proof. By direct calculation, we obtain
\[ \mathcal{P}(Z,W,s) = \sum_{g\in\Gamma_q} \bigl\{ \det(\Im(Z_1))^s \det(\Im(W))^s \left|\det(Z_1+W)\right|^{-2s} \det(Z_1+W)^{-k} \mathfrak{Q}(Q(Z,W),s) \bigr\}|_{\rho_{q*}} g. \]
Therefore
\begin{align*}
   & (f_*,\mathcal{P}(-\overline{Z},*,\overline{s})) \tag{5.4}
\\ & \qquad = \int_{\Gamma_q\backslash\mathfrak{H}_q} \langle \rho_{q*}(\sqrt{\Im(W)})f(W)_*, \rho_{q*}(\sqrt{\Im(W)})\mathcal{P}(-\overline{Z},W,\overline{s})\rangle \det(\Im(W))^{-q-1}\, dW
\\ & \qquad = 2 \int_{\mathfrak{H}_q} \det(\Im(Z_1))^s \det(\Im(W))^s \left|\det(-\overline{Z_1}+W)\right|^{-2s} \det(-\overline{Z_1}+W)^{-k}
\\ & \qquad\qquad \langle \rho_{q*}(\Im(W))f(W)_*, \mathfrak{Q}(Q(-\overline{Z},W),\overline{s})\rangle \det(\Im(W))^{-q-1}\, dW.
\end{align*}

 We compute the integral (5.4) according to Klingen [9, Section 1] (see also [1], [3], [17]). First we can choose $F=\begin{pmatrix} F_1^{(q)} & 0 \\ F_2 & F_3 \end{pmatrix}\in GL(p,\mathbb{R})$ satisfying $\Im(Z)=F^{-1}{}^t\!F^{-1}$. Then we have $\Im(Z_1)=F_1^{-1}{}^t\!F_1^{-1}$. Next we put
\[ S:=F_1(W-Z_1)(W-\overline{Z_1})F_1^{-1}. \]
Then we have
\[ \det(\Im(Z))^s \det(\Im(W))^s \left|\det(-\overline{Z_1}+W)\right|^{-2s}\det(\Im(W))^{-q-1}\,dW = 2^{q(q+1-2s)}\det(1_q-\overline{S}S)^{s-q-1}\,dS. \]
Furthermore we put
\[ \hat{f}(S):=\rho_q(W-\overline{Z_1})f(W). \]
Then the integral (5.4) is equal to
\[ 2^{q(q+1-2s)-2\left|\rho_q\right|+1} \int_{S_q} \langle \hat{f}(S)_*, \rho_{q*}({}^t\!F_1(1_q-\overline{S}S)F_1)\mathfrak{Q}(Q_2,\overline{s})\rangle \det(1_q-\overline{S}S)^{s-q-1}\, dS, \]
where
\[ Q_2:=-\frac{1}{2i}\begin{pmatrix} \begin{pmatrix} {}^t\!F_1SF_1 & 0 \\ 0 & 0 \end{pmatrix} + \begin{pmatrix} {}^t\!F_2F_2 & {}^t\!F_2F_3 \\ {}^t\!F_3F_2 & {}^t\!F_3F_3 \end{pmatrix} & -2i\begin{pmatrix} 1_q \\ 0 \end{pmatrix} \\ -2i\begin{pmatrix} 1_q & 0 \end{pmatrix} & 2^2F_1^{-1}\overline{S}(1_q-S\overline{S})^{-1}{}^t\!F_1^{-1} \end{pmatrix}. \]

 Next we consider the integral
\[ \int_{S_q} \langle \hat{f}(S)_*, \rho_{q*}({}^t\!F_1(1_q-\overline{S}S)F_1)\mathfrak{Q}(Q_2,\overline{s})\rangle \det(1_q-\overline{S}S)^{s-q-1}\, dS. \tag{5.5} \]
For a complex variable $t$ with $\left|t\right|\leq 1$, we put
\[ g(t):=\hat{f}(tS). \]
Then $g(t)$ has Taylor expansion
\[ g(t) = \sum_{\nu=0}^\infty \hat{f}_\nu(S) t^{\nu}. \]
Therefore we have
\[ \hat{f}_\nu(tS) = t^{\nu}\hat{f}_\nu(S) \]
and
\[ \hat{f}(S) = \sum_{\nu=0}^\infty \hat{f}_\nu(S). \]

On the other hand, we write $\mathfrak{Q}(X,s)$ as
\[ \mathfrak{Q}(X,s) = \sum_{\genfrac{}{}{0pt}{}{L\in\mathcal{L}}{\overline{L}=I}} c(L,s) X^L. \]
We put $Q_2 = -\frac{1}{2i}(R_{10}+R_{11}+R_2+{}^t\!R_2+R_3)$ with
\[ R_{10} := \begin{pmatrix} \begin{pmatrix} {}^t\!F_1SF_1 & 0 \\ 0 & 0 \end{pmatrix} & 0 \\ 0 & 0 \end{pmatrix}, \quad R_{11} := \begin{pmatrix} \begin{pmatrix} {}^t\!F_2F_2 & {}^t\!F_2F_3 \\ {}^t\!F_3F_2 & {}^t\!F_3F_3 \end{pmatrix} & 0 \\ 0 & 0 \end{pmatrix}, \]
\[ R_2 := \begin{pmatrix} 0 & -2i\begin{pmatrix} 1_q \\ 0 \end{pmatrix} \\ 0 & 0 \end{pmatrix} \quad\text{and}\quad R_3 := \begin{pmatrix} 0 & 0 \\ 0 & 2^2F_1^{-1}\overline{S}(1_q-S\overline{S})^{-1}{}^t\!F_1^{-1} \end{pmatrix}. \]
Then for $L$, there exist $L_1\in\mathcal{L}^{**}\sqcup\{\emptyset\}$, $L_2\in\mathcal{L}^*_*\sqcup\{\emptyset\}$ and $L_3\in\mathcal{L}_{**}\sqcup\{\emptyset\}$ such that
\[ Q_2^L = \left(-\frac{1}{2i}\right)^{\# L} (R_{10}+R_{11})^{L_1} R_2^{L_2} R_3^{L_3}. \]
Furthermore we can expand $(R_{10}+R_{11})^{L_1}$ as
\[ (R_{10}+R_{11})^{L_1} = \sum_{L_{10}\sqcup L_{11}=L_1} R_{10}^{L_{10}} R_{11}^{L_{11}}. \]
Therefore $\mathfrak{Q}(Q_2,\overline{s})$ can be expressed as
\[ \mathfrak{Q}(Q_2,\overline{s}) = \sum c(L,\overline{s}) \left(-\frac{1}{2i}\right)^{\# L} R_{10}^{L_{10}} R_{11}^{L_{11}} R_2^{L_2} R_3^{L_3}, \]
where the summation runs over
\[ \left\{ (L_{10},L_{11},L_2,L_3) \biggm| \genfrac{}{}{0pt}{}{L=L_{10}\sqcup L_{11}\sqcup L_2\sqcup L_3\in\mathcal{L},\quad\overline{L}=I,}{L_{10}, L_{11}\in\mathcal{L}^{**}\sqcup\{\emptyset\},\quad L_2\in\mathcal{L}^*_*\sqcup\{\emptyset\},\quad L_3\in\mathcal{L}_{**}\sqcup\{\emptyset\}} \right\}. \]

For $\nu$ and $L_{10}$, $L_{11}$, $L_2$, $L_3$, we consider the integral
\[ \int_{S_q} \langle \hat{f}_\nu(S)_*, \rho_{q*}({}^t\!F_1(1_q-\overline{S}S)F_1) R_{10}^{L_{10}} R_{11}^{L_{11}} R_2^{L_2} R_3^{L_3}\rangle \det(1_q-\overline{S}S)^{s-q-1}\, dS. \tag{5.6} \]
Substituting $e^{i\psi}S$ with a real variable $\psi$ for $S$, (5.6) is equal to
\[ e^{i\psi(\nu+\# L_3-\# L_{10})} \int_{S_q} \langle \hat{f}_\nu(S)_*, \rho_{q*}({}^t\!F_1(1_q-\overline{S}S)F_1) R_{10}^{L_{10}} R_{11}^{L_{11}} R_2^{L_2} R_3^{L_3}\rangle \det(1_q-\overline{S}S)^{s-q-1}\, dS. \]
Therefore if $\nu+\# L_3-\# L_{10}\neq 0$ then the integral vanishes. Since $\# L_{10}+\# L_{11} =\#L_3$, $\# L_3-\# L_{10}\geq 0$. Then the integral vanishes unless $\nu=0$ and $L_{11}=\emptyset$.

Consequently, the integral (5.5) is equal to
\[ \int_{S_q} \langle \hat{f_0}(S)_*, \rho_{q*}({}^t\!F_1(1_q-\overline{S}S)F_1)\mathfrak{Q}(Q_3,\overline{s})\rangle \det(1_q-\overline{S}S)^{s-q-1}\, dS, \tag{5.7} \]
where
\[ Q_3:=-\frac{1}{2i}\begin{pmatrix} \begin{pmatrix} {}^t\!F_1SF_1 & 0 \\ 0 & 0 \end{pmatrix} & -2i\begin{pmatrix} 1_q \\ 0 \end{pmatrix} \\ -2i\begin{pmatrix} 1_q & 0 \end{pmatrix} & 2^2F_1^{-1}\overline{S}(1_q-S\overline{S})^{-1}{}^t\!F_1^{-1} \end{pmatrix}. \]
Since $\hat{f}_0(S)=\hat{f}(0)=\rho_q(Z_1-\overline{Z_1})f(Z_1)$ and
\[ \mathfrak{Q}(Q_3,\overline{s}) = \mathfrak{Q}(Q_4,\overline{s}) \quad\text{with}\quad Q_4:=-\frac{1}{2i}\begin{pmatrix} {}^t\!F_1SF_1 & -2i\,1_q \\ -2i\,1_q & 2^2F_1^{-1}\overline{S}(1_q-S\overline{S})^{-1}{}^t\!F_1^{-1} \end{pmatrix}, \]
(5.7) is equal to
\begin{align*}
    & \int_{S_q} \langle \rho_{q*}(Z_1-\overline{Z_1})f(Z_1)_*, \rho_{q*}({}^t\!F_1(1_q-\overline{S}S)F_1)\mathfrak{Q}(Q_4,\overline{s})\rangle \det(1_q-\overline{S}S)^{s-q-1}\, dS
\\  & = (2i)^{\left|\rho_q\right|} \int_{S_q} \langle \rho_{q*}((1_q-\overline{S}S){}^t\!F_1^{-1})f(Z_1)_*, \rho_{q*}(F_1)\mathfrak{Q}(Q_4,\overline{s})\rangle \det(1_q-\overline{S}S)^{s-q-1}\, dS
\\  & = (2i)^{\left|\rho_q\right|} c(s,\rho_q) f(Z_1)^*.
\end{align*}
Thus Lemma 5.3 is proved. \qed

\bigskip

Summing up Proposition 5.1, (5.2), (5.3) and Lemma 5.3, we have
\begin{align*}
    & \left( f_*,(\mathfrak{P}(\partial)E^{p+q}_k)\Bigl(\begin{pmatrix} -\overline{Z}^{(p)} & 0 \\ 0 & * \end{pmatrix},\overline{s} \Bigr) \right)
\\  & = \sum_{T\in\mathbb{T}^{(q)}} (f_*, P_q(-\overline{Z},*,T,\overline{s}))
\\  & = \sum_{T\in\mathbb{T}^{(q)}} \lambda(T,f) \det(T)^{-k-2s}
\\  & \qquad\qquad \cdot 2^{q(q+1-2s)+1}\, (2^{-1}i)^{\left|\rho_q\right|}\, c(s,\rho_q) \sum_{g_1''\in\Gamma_{p,q}\backslash\Gamma_p} \left(\frac{\det(\Im(Z))}{\det(\Im(Z_1))}\right)^{s} f(Z_1)^* \Bigm|_{\rho_p^*} g_1''
\\  & = 2^{q(q+1-2s)+1}\, (2^{-1}i)^{\left|\rho_q\right|}\, c(s,\rho_q)\, D(k+2s,f)\, [f]^{\rho_p}_{\rho_q}(Z,s)^*.
\end{align*}
Thus Theorem 5.2 is proved. \qed

\bigskip

 For a cuspidal eigenform $f\in\mathfrak{S}_{\rho_q}$, we put
\[ c(s,f) := 2^{q(q+1-2s)+1}\, (2^{-1}i)^{\left|\rho_q\right|}\, c(s,\rho_q)\, D(k+2s,f). \tag{5.8} \]
Then from Theorem 2.1, Theorem 2.2 and Theorem 5.2, we have the following:

\bigskip
{\bf Proposition 5.4.}
{\sl
 For a cuspidal eigenform\/ $f\in\mathfrak{S}_{\rho_q}$, $c(s,f)[f]^{\rho_p}_{\rho_q}(Z,s)$ has meromorphic continuation to the whole\/ $s$-plane.
}

\bigskip
{\bf Corollary 5.5.}
{\sl Let\/ $f\in\mathfrak{S}_{\rho_q}$ be a cuspidal eigenform. We assume that\/ $c(s,\rho_q)$ is a meromorphic function in\/ $s$ and it is not identically zero. Then\/ $D(2s+k,f)[f]^{\rho_p}_{\rho_q}(Z,s)$ (or\/ $L(2s+k-q,f,\underline{\rm St})[f]^{\rho_p}_{\rho_q}(Z,s)$) has meromorphic continuation to the whole\/ $s$-plane.
}

\bigskip
{\bf Proposition 5.6.}
{\sl
 For a cuspidal eigenform\/ $f\in\mathfrak{S}_{\rho_q}$, $\gamma^{p+q}_k(s) c\bigl(s-\frac{k}{2},f\bigr)[f]^{\rho_p}_{\rho_q}\bigl(Z,s-\frac{k}{2}\bigr)$ is an entire function in\/ $s$, and it is invariant under the substitution $s\mapsto\frac{p+q-1}{2}-s$.
}

\bigskip
\noindent{\bf 6. Analytic properties of standard $L$-functions and Klingen Eisenstein series}

In this section we investigate the standard $L$-function $L(s,f,\underline{\rm St})$ and Klingen Eisenstein series $[f]^{\rho_p}_{\rho_q}(Z,s)$ for a cuspidal eigenform $f\in\mathfrak{S}_q$. First from (2.3), we have
\[ L(s,f,\underline{\rm St}) = \zeta(s+q)\prod_{j=1}^q\zeta(2s+2q-2j)D(s+q,f). \]
Furthermore we define the completed standard $L$-function $\Lambda(s,f,\underline{\rm St})$ by
\[ \Lambda(s,f,\underline{\rm St}):=\Gamma_{\rho_q}(s)L(s,f,\underline{\rm St}), \]
where
\[ \Gamma_{\rho_q}(s):=\Gamma_{\mathbb{R}}(s+\varepsilon_q)\prod_{j=1}^q\Gamma_{\mathbb{C}}(s+k+\lambda_j-j), \]
\[ \Gamma_{\mathbb{R}}(s):=\pi^{-s/2}\Gamma\left(\frac{s}{2}\right),\quad \Gamma_{\mathbb{C}}(s):=2(2\pi)^{-s}\Gamma(s), \]
and
\[ \varepsilon_q := 
    \begin{cases} 0 & \quad \text{if\/ $q$ even,} \\ 1 & \quad \text{if\/ $q$ odd.} \end{cases} \]
Then by direct calculation, we obtain
\begin{align*}
    \Lambda(s,f,\underline{\rm St}) &= 2^{q(2-s+k)-\left|\rho_q\right|} \pi^{(2q^2+q-\varepsilon_q)/2 - \left|\rho_q\right|} \prod_{j=1}^q \frac{\Gamma(s+k+\lambda_j-j)}{\Gamma(s+q+k+1-2j)}
\\                                  &\quad \cdot \gamma_{q,q}(s)^{-1}\, \gamma^{2q}_k\bigl(\frac{s+q}{2}\bigr) D(s+q,f),
\end{align*}
where
\[ \gamma_{p,q}(s) := \begin{cases} \frac{\displaystyle{\Gamma_p\Bigl(\frac{s+q}{2}\Bigr)}}{\displaystyle{\Gamma_p\Bigl(\frac{s}{2}\Bigr)}} & \quad\text{if $q$ even}, \\ \frac{\displaystyle{\Gamma_{p-1}\Bigl(\frac{s+q}{2}\Bigr)}}{\displaystyle{\Gamma_{p-1}\Bigl(\frac{s-1}{2}\Bigr)}} & \quad\text{if $q$ odd.} \end{cases} \]
We note that $\gamma_{p,q}(s)$ is a polynomial in $s$ and
\[ \deg\gamma_{p,q} = \begin{cases} \displaystyle{\frac{pq}{2}} & \quad\text{if $q$ even}, \\ \displaystyle{\frac{(p-1)(q+1)}{2}} & \quad\text{if $q$ odd}. \end{cases} \]
Furthermore $\gamma_{p,q}(s)$ satisfies the functional equation
\[ \gamma_{p,q}(s) = (-1)^{\deg \gamma_{p,q}} \, \gamma_{p,q}(p-q+1-s). \]

Next we define a modified Klingen Eisenstein series $\mathcal{E}^{\rho_p}_{\rho_q}(Z,s,f)$ by
\[ \mathcal{E}^{\rho_p}_{\rho_q}(Z,s,f) := \Lambda(s,f,\underline{\rm St}) \cdot \frac{\displaystyle{\Gamma_{p-q}\Bigl(\frac{s-q+k}{2}\Bigr)}}{\displaystyle{\Gamma_{p-q}\Bigl(\frac{s+\varepsilon_q}{2}\Bigr)}} \cdot \prod_{j=q+1}^{[(p+q)/2]} \xi(2s+2q-2j) \cdot [f]^{\rho_p}_{\rho_q}\bigl(Z,\frac{s+q-k}{2}\bigr). \]
From Theorem 5.2 and the definition of the Eisenstein series $\mathbb{E}^{p+q}_k(\mathfrak{Z},s)$, we obtain
\begin{align*}
    \mathcal{E}^{\rho_p}_{\rho_q}(Z,s,f)^* &= 2^{q-1} \pi^{(2q^2+q-\varepsilon_q)/2} (\pi i)^{-\left|\rho_q\right|} \gamma_{p,q}(s)^{-1} \tag{6.1}
\\  & \qquad \cdot \prod_{j=1}^q \frac{\Gamma(s+k+\lambda_j-j)}{\Gamma(s+q+k+1-2j)} \cdot c\bigl(\frac{s+q-k}{2},\rho_q\bigr)^{-1}
\\  & \qquad \cdot \left( f_*,(\mathfrak{P}(\partial)\mathbb{E}^{p+q}_k)\Bigl(\begin{pmatrix} -\overline{Z} & 0 \\ 0 & * \end{pmatrix},\frac{\overline{s}+q}{2} \Bigr) \right).
\end{align*}
We expect that
\[ \prod_{j=1}^q \frac{\Gamma(s+k+\lambda_j-j)}{\Gamma(s+q+k+1-2j)} \cdot c\bigl(\frac{s+q-k}{2},\rho_q\bigr)^{-1} \]
does not depend on $s$, i.e.,

\bigskip
{\bf Conjecture 6.1.}
{\sl
There exists a non-zero constant\/ $c_{\rho_q,\mathfrak{P}}$ depending only on\/ $\rho_q$ and\/ $\mathfrak{P}$ such that
\[ c\bigl(\frac{s+q-k}{2},\rho_q) = c_{\rho_q,\mathfrak{P}} \cdot \prod_{j=1}^q \frac{\Gamma(s+k+\lambda_j-j)}{\Gamma(s+q+k+1-2j)}. \]
}

\noindent Remarks. (1) In symmetric or alternating tensor valued case, for a suitable $\mathfrak{P}$, this conjecture is true.

(2) If this conjecture is true, by Corollary 5.5, $\mathcal{E}^{\rho_p}_{\rho_q}(Z,s,f)$ has meromorphic continuation to the whole $s$-plane.

\bigskip

If Conjecture 6.1 is true, (6.1) can be written as
\begin{align*}
    \mathcal{E}^{\rho_p}_{\rho_q}(Z,s,f)^* &= 2^{q-1} \pi^{(2q^2+q-\varepsilon_q)/2} (\pi i)^{-\left|\rho_q\right|} \gamma_{p,q}(s)^{-1} c_{\rho_q,\mathfrak{P}}^{-1} \tag{6.2}
\\  & \qquad \cdot \left( f_*,(\mathfrak{P}(\partial)\mathbb{E}^{p+q}_k)\Bigl(\begin{pmatrix} -\overline{Z} & 0 \\ 0 & * \end{pmatrix},\frac{\overline{s}+q}{2} \Bigr) \right).
\end{align*}
Using the functional equation of $\mathbb{E}^{p+q}_k(\mathfrak{Z},s)$ and $\gamma_{p,q}(s)$, we have the following:

\bigskip
{\bf Proposition 6.2.}
{\sl
Let\/ $f\in\mathfrak{S}_q$ be a cuspidal eigenform. If Conjecture 6.1 is true, then\/ $\mathcal{E}^{\rho_p}_{\rho_q}(Z,s,f)$ satisfies the functional equation
\[ \mathcal{E}^{\rho_p}_{\rho_q}(Z,s,f) = (-1)^{\deg\gamma_{p,q}}\, \mathcal{E}^{\rho_p}_{\rho_q}(Z,p-q+1-s,f). \]
In particular, the completed standard\/ $L$-function\/ $\Lambda(s,f,\underline{\rm St})$ satisfies the functional equation
\[ \Lambda(s,f,\underline{\rm St}) = \Lambda(1-s,f,\underline{\rm St}). \]
}

Next we investigate the location of poles of $\mathcal{E}^{\rho_p}_{\rho_q}(Z,s,f)$. The result is:

\bigskip
{\bf Theorem 6.3.}
{\sl
Let\/ $f\in\mathfrak{S}_q$ be a cuspidal eigenform. We assume that Conjecture 6.1 is true.

{\rm (i)} If\/ $k\geq\displaystyle{\frac{p+q}{2}}$ and\/ $p+q\not\equiv 0$ $(\mod 4)$, then\/ $\mathcal{E}^{\rho_p}_{\rho_q}(Z,s,f)$ is entire.

{\rm (ii)} If\/ $k\geq\displaystyle{\frac{p+q}{2}}$ and\/ $p+q\equiv 0$ $(\mod 4)$, then\/ $\mathcal{E}^{\rho_p}_{\rho_q}(Z,s,f)$ is entire with the exception of the following cases:
\[ \begin{matrix}
    \text{\rm (1)}\quad p-q=0 & \text{and\/ $q$ is even},
\\  \text{\rm (2)}\quad p-q=2 & \text{and\/ $q$ is odd}.
   \end{matrix} \]
In these cases,\/ $\mathcal{E}^{\rho_p}_{\rho_q}(Z,s,f)$ is holomorphic except for possible simple poles at\/ $s=\displaystyle{\frac{p-q}{2}}$ and\/ $s=1+\displaystyle{\frac{p-q}{2}}$.

{\rm (iii)} If\/ $k<\displaystyle{\frac{p+q}{2}}$ and\/ $k\geq q+\varepsilon_q+2$, then\/ $\mathcal{E}^{\rho_p}_{\rho_q}(Z,s,f)$ is entire.

{\rm (iv)} If\/ $k<\displaystyle{\frac{p+q}{2}}$ and\/ $k\leq q+\varepsilon_q$, then\/ $\mathcal{E}^{\rho_p}_{\rho_q}(Z,s,f)$ is holomorphic except for possible poles at
\[ k-q\leq s\leq p-k+1, \quad s\in\mathbb{Z}. \]
Furthermore for\/ $j=0$, $1$, $\ldots$, $[\frac{p+q}{2}]-k$,
\begin{align*}
    {\rm ord}\bigl(k-q+j;\mathcal{E}^{\rho_p}_{\rho_q}(Z,*,f)\bigr) &\leq \min\left\{\left[\frac{j}{2}\right], \frac{q+\varepsilon_q-k}{2}\right\}+1,
\\  {\rm ord}\bigl(p-k+1-j;\mathcal{E}^{\rho_p}_{\rho_q}(Z,*,f)\bigr) &\leq \min\left\{\left[\frac{j}{2}\right], \frac{q+\varepsilon_q-k}{2}\right\}+1.
\end{align*}
If\/ $p+q$ is odd then
\[ {\rm ord}\Bigl(\frac{p-q+1}{2};\mathcal{E}^{\rho_p}_{\rho_q}(Z,*,f)\Bigr) \leq \frac{q+\varepsilon_q-k}{2}. \]
Here ${\rm ord}(s;\phi)$ denotes the order of poles at\/ $s$ for a meromorphic function\/ $\phi$.
}

\bigskip
{\bf Corollary 6.4.}
{\sl
Let\/ $f\in\mathfrak{S}_q$ be a cuspidal eigenform. We assume that Conjecture 6.1 is true.

{\rm (i)} If\/ $k\geq q$, then\/ $\Lambda(s,f,\underline{\rm St})$ is holomorphic except for possible simple poles at\/ $s=0$ and\/ $s=1$. Furthermore if $q$ is odd then $\Lambda(s,f,\underline{\rm St})$ is entire.

{\rm (ii)} If\/ $k<q$, then\/ $\Lambda(s,f,\underline{\rm St})$ is holomorphic except for possible poles at
\[ k-q\leq s\leq q-k+1, \quad s\in\mathbb{Z}. \]
Furthermore for\/ $j=0$, $1$, $\ldots$, $q-k$,
\begin{align*}
    {\rm ord}\bigl(k-q+j;\Lambda(*,f,\underline{\rm St})\bigr) &\leq \left[\frac{j}{2}\right]+1,
\\  {\rm ord}\bigl(q-k+1-j;\Lambda(*,f,\underline{\rm St})\bigr) &\leq \left[\frac{j}{2}\right]+1.
\end{align*}
}

\noindent Remark. By Weissauer [20], $\Lambda(s,f,\underline{\rm St})$ has a pole at $s=1$ (or equivalently, at $s=0$) if and only if $f$ can be written as a linear combination of some theta series attached to positive definite even unimodular matrices of size $2q$ and pluri-harmonic polynomials (see [20], [15], [17]). Therefore if $k\geq q$ and $q\not\equiv 0$ $(\mod 4)$, then $\Lambda(s,f,\underline{\rm St})$ is entire.

\bigskip

Proof of Theorem 6.3. Theorem 6.3 is proved in the same way as that by Mizumoto [15, Theorem 3]. \qed

\bigskip
\noindent{\bf 7. Algebraicity results for Siegel modular forms and standard $L$-functions}

Each $f\in\mathfrak{M}_{\rho_q}$ has a Fourier expansion
\[ f(Z) = \sum_{R} a(R,f) e^{2\pi i{\rm tr}(RZ)}, \]
where $R$ runs over all symmetric positive semi-definite semi-integral matrices of size $q$. Let $K$ be a subfield of $\mathbb{C}$. We put
\[ \mathfrak{M}_{\rho_q}(K) := \{ f\in\mathfrak{M}_{\rho_q} \mid a(R,f)\in V_q(K) \quad\text{for all $R$} \}, \]
where $V_q(K)$ denotes the set of $K$-rational points of $V_q$. For any $\mathbb{C}$-subspace $\mathfrak{W}$ of $\mathfrak{M}_{\rho_q}$, put
\[ \mathfrak{W}(K) := \mathfrak{W} \cap \mathfrak{M}_{\rho_q}(K). \]

Let $\mathcal{H}_K:=\mathcal{H}_K(\Gamma_q,G^+Sp(q,\mathbb{Q}))$ be the Hecke algebra. Let
\[ t_{\rho_q}\: \mathcal{H}_K \to {\rm End}_{\mathbb{C}}(\mathfrak{S}_{\rho_q}) \]
be the usual representation. We put
\[ \mathbb{T}_K := t_{\rho_q}(\mathcal{H}_K). \]
For any $\lambda\in\hat{\mathbb{T}}_{\mathbb{C}}:={\rm Hom}_{\mathbb{C}\text{-\rm alg}}(\mathbb{T}_{\mathbb{C}},\mathbb{C})$, put
\[ \mathfrak{S}_{\rho_q,\lambda} := \{f\in\mathfrak{S}_{\rho_q} \mid Tf=\lambda(T)f \quad\text{for all $T\in\mathbb{T}_{\mathbb{C}}$} \} \]
and define the field
\[ \mathbb{Q}(\lambda):=\mathbb{Q}(\lambda(T) \mid T\in\mathbb{T}_{\mathbb{Q}}). \]
Let
\[ \Lambda := \{ \lambda\in\hat{\mathbb{T}}_{\mathbb{C}} \mid \mathfrak{S}_{\rho_q,\lambda}\neq \{0\} \}. \]
Then
\[ \mathfrak{S}_{\rho_q} = \bigoplus_{\lambda\in\Lambda} \mathfrak{S}_{\rho_q,\lambda}. \]

Let $f\in\mathfrak{S}_{\rho_q}$ be an eigenform. Then there exists $\lambda\in\Lambda$ such that $f\in\mathfrak{S}_{\rho_q,\lambda}$. Put $\mathbb{Q}(f):=\mathbb{Q}(\lambda)$. Note that $\mathbb{Q}(f)$ is a subfield of $\mathbb{R}$ since $\lambda(T)$ is hermitian for $T\in\mathbb{T}_{\mathbb{Q}}$.

Using the same method as in Mizumoto [15], we obtain the following theorems: (see also [19], [11])

\bigskip
{\bf Theorem 7.1}.
{\sl
We assume that\/ $k>q+1$ or\/ $k\equiv 0$ $(\mod 4)$.

(1) Let
\[ \mathfrak{V}_{\rho_q} := \bigoplus_{\genfrac{}{}{0pt}{}{\lambda\in\Lambda}{c(\lambda)\neq 0}} \mathfrak{S}_{\rho_q,\lambda}. \]
Then\/ $\mathfrak{V}_{\rho_q}$ has a basis consisting of elements of\/ $\mathfrak{S}_{\rho_q}(\mathbb{Q})$. In particular,\/ ${\rm Aut}(\mathbb{C})$ acts on\/ $\mathfrak{V}_{\rho_q}$ by\/ $f\mapsto f^{\sigma}$.

(2) Let\/ $\lambda\in\Lambda$ with $c(\lambda)\neq 0$. Then

\noindent{\rm (i)}\/ $\mathbb{Q}(\lambda)$ is a totally real finite extension of\/ $\mathbb{Q}$ with
\[ [\mathbb{Q}(\lambda):\mathbb{Q}] \leq \#\Lambda. \]

\noindent{\rm (ii)}\/ $\mathfrak{S}_{\rho_q,\lambda}$ has a basis consisting of elements of $\mathfrak{S}_{\rho_q,\lambda}(\mathbb{Q}(\lambda))$. 

\noindent{\rm (iii)} For\/ $0\neq f\in\mathfrak{S}_{\rho_q,\lambda}(\mathbb{Q}(\lambda))$, there exists an orthogonal basis\/ $\{f_1,f_2,\ldots,f_d\}$ of $\mathfrak{S}_{\rho_q,\lambda}$ such that
\[ f_1=f \quad\text{and}\quad f_1,f_2,\ldots,f_d\in\mathfrak{S}_{\rho_q,\lambda}(\mathbb{Q}(\lambda)). \]

\noindent{\rm (iv)} For\/ $0\neq f\in\mathfrak{S}_{\rho_q,\lambda}(\mathbb{Q}(\lambda))$ and\/ $\sigma\in{\rm Aut}(\mathbb{C})$,
\[ \left(\frac{c(f)}{(f,f)}\right)^{\sigma} = \frac{c(f^\sigma)}{(f^\sigma,f^\sigma)}. \]
}

\bigskip
{\bf Theorem 7.2}.
{\sl
Let\/ $k>(p+q+3)/2$ or\/ $k\equiv 0$ $(\mod 4)$. Let\/ $\lambda\in\Lambda$ with\/ $c(\lambda)\neq 0$. Then for $0\neq f\in\mathfrak{S}_{\rho_q,\lambda}(\mathbb{Q}(\lambda))$ and $\sigma\in{\rm Aut}(\mathbb{C})$,
\[ ([f]^{\rho_p}_{\rho_q})^\sigma = [f^\sigma]^{\rho_p}_{\rho_q}. \]
}

\bigskip
{\bf Theorem 7.3}.
{\sl 
Let\/ $f\in\mathfrak{S}_{\rho_q}$ be an eigenform. Suppose that each Fourier coefficient of\/ $f$ belongs to $V_q(\mathbb{Q}(f))$. Let\/ $s\in\mathbb{Z}$ be such that
\[ 1\leq s\leq k-q \quad\text{and}\quad s\equiv q \enskip (\mod 2). \]
We assume
\[ \text{$q\equiv 3$ $(\mod 4)$ \quad if\/ $s=1$}. \]
Let
\[ A(f) := \frac{c(0,\rho_q) \cdot L(s,f,\underline{\rm St})}{\pi^{\left|\rho_q\right|+s(q+1)}(f,f)}. \]
Then we have
\[ A(f)^{\sigma} = A(f^\sigma) \quad\text{for all\/ $\sigma\in{\rm Aut}(\mathbb{C})$}. \]
In particular,
\[ A(f)\in\mathbb{Q}(f). \]
}

\bigskip
\noindent{\bf References} \par

{
\parindent=2pc

\begin{itemize}

\item[{[1]}] S.~B\"ocherer, \"Uber die Funktionalgleichung automorpher $L$-Funktionen zur Siegelschen Modulgruppe, J. Reine Angew. Math., {\bf 362} (1985), 146--168.

\item[{[2]}] S.~B\"ocherer, Ein Rationalit\"atssatz f\"ur formale Heckereihen zur Siegelschen Modulgruppe, Abh. Math. Sem. Univ. Hamburg, {\bf 56} (1986), 35--47.

\item[{[3]}] S.~B\"ocherer, T.~Satoh, and T.~Yamazaki, On the pullback of a differential operator and its application to vector valued Eisenstein series, Comment. Math. Univ. St. Pauli, {\bf 42} (1992), 1--22.

\item[{[4]}] B.~Diehl, Die analytische Fortsetzung der Eisensteinreihe zur Siegelschen Modulgruppe, J. Reine Angew. Math., {\bf 317} (1980), 40--73.

\item[{[5]}] P.~Feit, Poles and residues of Eisenstein series for symplectic and unitary groups, Mem. Amer. Math. Soc. {\bf 61} (1986), no. 346.

\item[{[6]}] P.~B.~Garrett, Pullbacks of Eisenstein series; applications, Progress in Math., {\bf 46} (1984), 114--137.

\item[{[7]}] T.~Ibukiyama, On differential operators on automorphic forms and invariant pluri-harmonic polynomials, Comment. Math. Univ. St. Pauli, {\bf 48} (1999), 103--118.

\item[{[8]}] V.~L.~Kalinin, Eisenstein series on the symplectic group, Math. USSR-Sb., {\bf 32} (1977), 449--476; English translation.

\item[{[9]}] H.~Klingen, \"Uber Poincar\'esche Reihen zur Siegelschen Modulgruppe, Math. Ann., {\bf 168} (1967), 157--170.

\item[{[10]}] H.~Klingen, Zum Darstellungssatz f\"ur Siegelschen Modulformen, Math. Z., {\bf 102} (1967), 30--43.

\item[{[11]}] N.~Kozima, On special values of standard $L$-functions attached to vector valued Siegel modular forms, Kodai Math. J., {\bf 23} (2000), 255--265.

\item[{[12]}] N.~Kozima, Standard $L$-functions attached to alternating tensor valued Siegel modular forms, Osaka J. Math., {\bf 39} (2002), 245--258.

\item[{[13]}] N.~Kozima, Garrett's pullback formula for vector valued Siegel modular forms, J. Number Theory, {\bf 128} (2008), 235--250.

\item[{[14]}] R.~P.~Langlands, On the functional equations satisfied by Eisenstein series, Lecture Notes in Math., {\bf 544}, Springer, Berlin Heidelberg New York, 1976.

\item[{[15]}] S.~Mizumoto, Poles and residues of standard $L$-functions attached to Siegel modular forms, Math. Ann., {\bf 289} (1991), 589--612.

\item[{[16]}] S.~Mizumoto, Eisenstein series for Siegel modular groups, Math. Ann., {\bf 297} (1993), 581--625.

\item[{[17]}] H.~Takayanagi, Vector valued Siegel modular forms and their $L$-functions; Application of a differential operator, Japan J. Math., {\bf 19} (1994), 251--297.

\item[{[18]}] H.~Takayanagi, On standard $L$-functions attached to ${\rm alt}^{n-1}(\mathbb{C}^n)$-valued Siegel modular forms, Osaka J. Math., {\bf 32} (1995), 547--563.

\item[{[19]}] Y.~Takei, On algebraicity of vector valued Siegel modular forms, Kodai Math. J., {\bf 15} (1992), 445--457.

\item[{[20]}] R.~Weissauer, Stabile Modulformen und Eisensteinreihen, Lecture Notes in Math., {\bf 1219}, Springer, Berlin Heidelberg New York, 1986.

\item[{[21]}] H.~Weyl, The Classical Groups, Princeton Univ. Press, 1946.

\end{itemize}

}

\end{document}